\documentclass[12pt]{article}
\usepackage{mathrsfs}
\usepackage{amsfonts}
= 0.0in \topmargin = 0.3in \oddsidemargin=0.1in
\def\Dj{\hbox{D\kern-.73em\raise.30ex\hbox{-}
\raise-.30ex\hbox{}}}
\def\dj{\hbox{d\kern-.33em\raise.80ex\hbox{-}
\raise-.80ex\hbox{\kern-.40em}}}
\usepackage{epsfig}
\usepackage{amsmath,amsthm,amsfonts,amssymb,amscd}
\usepackage{graphicx,ifpdf}

\begin{document}
\baselineskip=0.30in

\newtheorem{lem}{Lemma}[section]
\newtheorem{thm}[lem]{Theorem}
\newtheorem{cor}[lem]{Corollary}
\newtheorem*{prop}{Proposition}
\newtheorem{con}[lem]{Conjecture}
\newtheorem{rem}[lem]{Remark}
\newtheorem{defi}[lem]{Definition}
\renewcommand\baselinestretch{1.2}
\def\pf{\noindent {\it Proof.} }
\def\qed{\hfill \rule{4pt}{7pt}}

%%%%%%\rule{0cm}{3.5cm}

\renewcommand\thefootnote{\fnsymbol{footnote}}

\begin{center} {\Large \bf On the minimal energy of tetracyclic graphs}
 \end{center}

\baselineskip=0.20in

%%%\vspace{2mm}

\baselineskip=0.20in \baselineskip=0.20in
\begin{center}
{ \small   Hongping Ma, Yongqiang
Bai\footnote{Corresponding author. %%%\indent
E-mail addresses: hpma@163.com (H. Ma), bmbai@163.com (Y. Bai)} \\[5pt]
School of Mathematics and Statistics, Jiangsu Normal University,\\ Xuzhou 221116, China\\
 }
\end{center}

\begin{abstract}
 The energy of a graph is defined as the sum of the absolute values of the eigenvalues of its adjacency
 matrix. In this paper, we characterize the tetracyclic graph of
 order $n$ with minimal energy. By this, the validity of a
 conjecture for the case $e=n+3$ proposed by Caporossi et al. \cite{CCGH} has been confirmed.
\\[2mm]
\end{abstract}

\baselineskip=0.245in

\section{Introduction}

 Let $G$ be a simple graph with $n$ vertices and $A(G)$ the adjacency matrix of
 $G$. The eigenvalues $\lambda_{1}, \lambda_{2},\ldots, \lambda_{n}$ of $A(G)$  are said to be the eigenvalues of the graph $G$.
The energy of $G$ is defined as
$$E=E(G)=\sum_{i=1}^{n}|\lambda_{i}|.$$
The characteristic polynomial of $A(G)$ is also called the
characteristic polynomial of $G$, denoted by
 $\phi(G,x)=\mbox{det}(xI-A(G))=\sum_{i=0}^{k}a_i(G)x^{n-i}$. Using
 these coefficients of $\phi(G,x)$, the energy of $G$ can be expressed as the Coulson integral formula
 \cite{GP}:
  \begin{eqnarray}
 E(G)=\frac{1}{2\pi}{\large\int}_{-\infty}^{+\infty}\frac{1}{x^{2}}\ln\left[\left(\sum\limits_{i=0}^{\lfloor
\frac{n}{2}\rfloor}
(-1)^ia_{2i}(G)x^{2i}\right)^2+\left(\sum\limits_{i=0}^{\lfloor
\frac{n}{2}\rfloor} (-1)^ia_{2i+1}(G)x^{2i+1}\right)^2\right]dx.
\label{energy-1}
\end{eqnarray}
For convenience, write $b_{2i}(G)=(-1)^ia_{2i}(G)$ and
$b_{2i+1}(G)=(-1)^ia_{2i+1}(G)$ for $0\leq i\leq\lfloor
\frac{n}{2}\rfloor$.

Since the energy of a graph can be used to approximate the total
$\pi$-electron energy of the molecular, it has been intensively
studied. For details on graph energy, we refer to the recent book
\cite{LSG} and reviews \cite{G2,GLZ}.

 One of the fundamental question that is encountered in the study of graph energy is which graphs (from a given class) have minimal and maximal energies.
 A large of number of papers were published on such extremal problems, see Chapter 7 in \cite{LSG}.

 A connected graph on $n$ vertices with $e$ edges is called an $(n,e)$-graph. We call an $(n,e)$-graph a unicyclic graph, a bicyclic graph,
 a tricyclic graph, and a tetracyclic graph if $e=n, n+1, n+2$ and $n+3$, respectively.
 Follow \cite{ZK}, let $S_{n,e}$ be the graph obtained by the star
 $S_n$ with $e-n+1$ additional edges all connected to the same vertex,
 and $B_{n,e}$ be the bipartite $(n,e)$-graph with two vertices on one side, one of which is connected to all vertices on the other side.

In \cite{CCGH}, Caporossi et al. gave the following conjecture:

\begin{con}\textnormal{\cite{CCGH}}\label{conjecture-minimal energy}
 Connected graphs $G$ with $n\geq 6$ vertices, $n-1\leq e\leq 2(n-2)$ edges and minimum energy are $S_{n,e}$ for $e\leq n+[(n-7)/2]$, and $B_{n,e}$ otherwise.
\end{con}

 This conjecture is true when $e=n-1$, $2(n-2)$ \cite{CCGH}%%%(Theorem1)
 , and when $e=n$ for $n\geq 6$ \cite{H1}. Li et al. \cite{LZW} showed that
 $B_{n,e}$ is the unique bipartite graph of order $n$ with minimal energy for
 $e\leq 2n-4$. Hou \cite{H2} proved that for $n\geq 6$, $B_{n,n+1}$ has the minimal
 energy among all bicyclic graphs of order $n$ with at most one odd
 cycle. Let $\mathcal{G}_{n,e}$ be the set of connected graphs with $n$ vertices and $e$ edges. Let $\mathcal{G}^1_{n,e}$ be the subset of $\mathcal{G}_{n,e}$
 which contains no disjoint two odd cycles of length $p$ and $q$ with $p+q\equiv 2$ (mod $4)$, and $\mathcal{G}^2_{n,e}=\mathcal {G}_{n,e}\setminus \mathcal{G}^1_{n,e}$.
 Zhang and Zhou \cite{ZZ} characterized the graphs with minimal, second-minimal and third-minimal energy in $\mathcal{G}^1_{n,n+1}$ for $n\geq 8$. Combining the results (Lemmas 5-9) in
\cite{ZZ} with
 the fact that $E(B_{n,n+1})<E(S_{n,n+1})$ for $5\leq n\leq 7$, we can
 deduce the following lemma.

\begin{lem}\textnormal{\cite{ZZ}}\label{bicyclic-minimal energy-1}
The graph with minimal energy in $\mathcal {G}^1_{n,n+1}$ is
$S_{n,n+1}$ for $n=4$ or $n\geq 8$, and $B_{n,n+1}$ for $5\leq n\leq
7$, respectively.
\end{lem}

 Li et al. \cite{LLZ} proved that $B_{n,n+2}$ has minimal energy in $\mathcal {G}^1_{n,n+2}$ for $7\leq n\leq 9$,
 %%%%and determined the graphs with minimal and second-minimal energy for $n\geq 10$ in the set $\mathcal {G}_{n,n+2}$ containing no four special graphs. Huo et al. solved this problem in \cite{HJL} and the results on minimal energy can be stated as follows.
 and for $n\geq 10$, they wanted to characterize the graphs with minimal and second-minimal energy in $\mathcal{G}^1_{n,n+2}$, but left four special graphs without determining
 their ordering. Huo et al. solved this problem in \cite{HJL}, and the results on minimal energy can be restated as follows.

\begin{lem}\label{tricyclic-minimal energy-1}
The graph with minimal energy in $\mathcal {G}^1_{n,n+2}$ is
$B_{n,n+2}$ for $7\leq n\leq 9$ \textnormal{\cite{LLZ}}, and
$S_{n,n+2}$ for $n\geq 10$ \textnormal{\cite{HJL}}, respectively.
\end{lem}

In \cite{ZK}, the authors claimed that they gave a complete solution
to conjecture \ref{conjecture-minimal energy} for $e=n+1$ and
$e=n+2$ by showing the following two results.

\begin{lem}\textnormal{(Theorem 1, \cite{ZK})}\label{bicyclic-minimal energy-2}
Let $G$ be a connected graph with $n$ vertices and $n+1$ edges. Then
$$E(G)\geq E(S_{n,n+1})$$ with equality if and only if $G\cong
S_{n,n+1}$.
\end{lem}

\begin{lem}\textnormal{(Theorem 2, \cite{ZK})}\label{tricyclic-minimal energy-2}
Let $G$ be a connected graph with $n$ vertices and $n+2$ edges. Then
$$E(G)\geq E(S_{n,n+2})$$ with equality if and only if $G\cong
S_{n,n+2}$.
\end{lem}

 Note that $E(B_{n,n+1})<E(S_{n,n+1})$ for $5\leq n\leq 7$, and $E(B_{n,n+2})<E(S_{n,n+2})$ for $6\leq n\leq
 9$. In addition, there is a little gap in the original proofs (even for large $n$) of Lemmas \ref{bicyclic-minimal energy-2} and \ref{tricyclic-minimal
 energy-2} in \cite{ZK}, respectively. For completeness, we will prove the following two results in Section 2.

 \begin{thm}\label{bicyclic-thm}
 $S_{n,n+1}$ if $n=4$ or $n\geq 8$, $B_{n,n+1}$ if $5\leq n\leq 7$ has minimal energy in $\mathcal{G}_{n,n+1}$.
\end{thm}

\begin{thm}\label{tricyclic-thm}
 The complete graph $K_4$ if $n=4$, $S_{n,n+2}$ if $n=5$ or $n\geq 10$, $B_{n,n+2}$ if $6\leq n\leq 9$ has minimal energy in $\mathcal{G}_{n,n+2}$.
 Furthermore, $S_{6,8}$ has second-minimal energy in $\mathcal{G}_{6,8}$.
\end{thm}

 Li and Li \cite{LL} discussed the graph with minimal energy in $\mathcal
 {G}^1_{n,n+3}$, and claimed that the graph with minimal energy in $\mathcal {G}^1_{n,n+3}$ is $B_{n,n+3}$ for $9\leq n\leq 17$, and $S_{n,n+3}$ for $n\geq 18$, respectively.
 Note that $E(S_{n,n+3})<E(B_{n,n+3})$ for $n\geq 12$. In Section 3, we will first illustrate the correct version of this
 result, and then we will show the following theorem.

\begin{thm}\label{tetracyclic-thm}
 The wheel graph $W_5$ if $n=5$, the complete bipartite graph $K_{3,3}$ if $n=6$, $B_{n,n+3}$ if $7\leq n\leq 11$, $S_{n,n+3}$ if $n\geq 12$ has minimal energy in $\mathcal{G}_{n,n+3}$.
 Furthermore, $S_{n,n+3}$ has second-minimal energy in $\mathcal{G}_{n,n+3}$ for $6\leq n\leq 7$.
\end{thm}

\begin{lem}\textnormal{\cite{ZK}}\label{lemma Sn,e and Bn,e}
 $E(S_{n,e}) < E(B_{n,e})$ if $n-1\leq e\leq \frac{3}{2}n-3$; $E(B_{n,e}) < E(S_{n,e})$ if $ \frac{3}{2}n-\frac{5}{2}\leq e\leq
 2n-4$.
\end{lem}

  From Lemma \ref{lemma Sn,e and Bn,e}, we know that the bound $e\leq n+[(n-7)/2]$ in Conjecture \ref{conjecture-minimal energy} should be understood that $e\leq n+ \lceil(n-7)/2\rceil$.
  With Theorems \ref{bicyclic-thm}, \ref{tricyclic-thm} and  \ref{tetracyclic-thm}, we give a complete solution to Conjecture \ref{conjecture-minimal energy} for $e=n+1,n+2$ and $n+3$.

\section{The graphs with minimal energy in $\mathcal{G}_{n,n+1}$ and
$\mathcal{G}_{n,n+2}$}

 The following three lemmas are need in the sequel.

\begin{lem}\textnormal{\cite{DS}}\label{edge-cut}
If $F$ is an edge cut of a simple graph $G$, then $E(G-F)\leq E(G)$,
where $G-F$ is the subgraph obtained from $G$ by deleting the edges
in $F$.
\end{lem}

\begin{lem}\textnormal{\cite{ZK}}\label{lemma in Zhang}
(1) Suppose that $n_1, n_2\geq 3$ and $n=n_1+n_2$. Then
$$E(S_{n_1,n_1 }\cup S_{n_2,n_2 })\geq E(S_{n-3,n-3}\cup C_3)$$ with
equality if and only if $\{n_1,n_2\}=\{3,n-3\}$.

(2) $E(S_{n-3,n-3}\cup C_3) > E(S_{n,n+1})$ for $n\geq 6$.

(3) $E(S_{n,n+1}) > E(S_{n,n})$ for $n\geq 4$.

(4) $E(S_{n-3,n-3}\cup C_3) > E(S_{n,n+2})$ for $n\geq 6$.
\end{lem}

\begin{lem}\label{unicyclic-minimal energy-1}
 (1) \textnormal{\cite{H1}} $S_{n,n}$ has minimal energy in $\mathcal{G}_{n,n}$ for $n=3$ or $n\geq 6$.

 (2) $B_{n,n}$ and $S_{n,n}$ have, respectively, minimal and second-minimal energy in $\mathcal{G}_{n,n}$ for $4\leq n\leq 5$.
 In particular, $S_{n,n}$ is the unique non-bipartite graph in $\mathcal{G}_{n,n}$ with minimal energy for $4\leq n\leq 5$.
\end{lem}

 \pf By Table 1 of \cite{CDS}, there are two $(4,4)$-graphs and five $(5,5)$-graphs. By simple computation, we can obtain the result (2).  \qed

\vskip 10pt

\noindent {\bf Proof of Theorem \ref{bicyclic-thm}: }
 By Lemma \ref{bicyclic-minimal energy-1}, it suffices to prove that
 $E(G)>E(S_{n,n+1})$ when $n=4$ or $n\geq 8$, and $E(G)>E(B_{n,n+1})$ when $5\leq n\leq 7$ for $G\in \mathcal{G}^2_{n,n+1}$.

 Suppose that $G\in \mathcal{G}^2_{n,n+1}$. As there is nothing to prove for the case $n\leq 5$, we suppose that $n\geq 6$. Then $G$ has a cut edge $f$ such that
 $G-f$ contains exactly two components, say $G_1$ and $G_2$, which are non-bipartite
 unicyclic graphs. Let $|V(G_1)|=n_1$, $|V(G_2)|=n_2$, and $n_1+n_2=n$.
 By Lemmas \ref{edge-cut}, \ref{lemma in Zhang} and \ref{unicyclic-minimal energy-1}, we have
\begin{eqnarray}
 E(G) & \geq  &  E(G_1\cup G_2) \label{Inq2}\\
      & \geq  &  E(S_{n_1,n_1}\cup S_{n_2,n_2 }) \label{Inq3}\\
      & \geq  &  E(S_{n-3,n-3}\cup C_3) \label{Inq4}\\
      &  >    & E(S_{n,n+1}) \label{Inq5}.
\end{eqnarray}
 %%%$$E(G)\geq E(G_1\cup G_2)\geq E(S_{n_1,n_1}\cup S_{n_2,n_2 })\geq E(S_{n-3,n-3}\cup C_3) > E(S_{n,n+1}).$$
 In particular, $E(G)> E(S_{n,n+1})> E(B_{n,n+1})$ for $6\leq n\leq
 7$.  The proof is thus complete. \qed

\begin{rem}\label{bicyclic-remark}
 The proof of Theorem \ref{bicyclic-thm} (for large $n$) is similar to that of Lemma \ref{bicyclic-minimal
 energy-2} except that in \cite{ZK}, the authors did not point out that $G_1$ and
 $G_2$ are non-bipartite unicyclic graphs. Without this assumption, we know that the
 inequality \eqref{Inq3} does not hold when $n_1$ or $n_2$ equals to $4$ or $5$ by Lemma \ref{unicyclic-minimal energy-1} (2).
 Moreover, the inequality $E(G_1\cup G_2)\geq E(S_{n-3,n-3}\cup C_3)$ does not hold. For
 example: $E(C_4\cup S_{n-4,n-4})< E(S_{n-3,n-3}\cup C_3)$ for
 $n\geq 7$, since $E(C_4)=E(C_3)=4$ and $E(S_{n-4,n-4})<E(S_{n-3,n-3})$ by Lemma \ref{edge-cut}.
\end{rem}

\begin{lem}\label{bicyclic-n-5,6,7}
 $S_{n,n+1}$ is the unique non-bipartite graph in $\mathcal{G}_{n,n+1}$ with minimal energy for $5\leq n\leq
 7$. Furthermore, $S_{n,n+1}$ has second-minimal energy in
 $\mathcal{G}_{n,n+1}$ for $n=5$ or $7$, and $S_{6,7}$ has third-minimal energy in
 $\mathcal{G}_{6,7}$.
\end{lem}

 \pf  By Table 1 of \cite{CDS}, there are five $(5,6)$-graphs. By simple calculation, we can prove the theorem for $n=5$.
 By Table 1 of \cite{CP}, there are 19 $(6,7)$-graphs. By direct computation, we can prove the theorem for $n=6$.
 By the results (Lemmas 5-9) in \cite{ZZ}, we can obtain that $S_{7,8}$ has second-minimal energy in
 $\mathcal{G}^1_{7,8}$. On the other hand, from the proof of
 Theorem \ref{bicyclic-thm}, $E(G)>E(S_{7,8})$ for $G\in\mathcal{G}^2_{7,8}$.
 Therefore $S_{7,8}$ has second-minimal energy in
 $\mathcal{G}_{7,8}$, and so the theorem is true for $n=7$.
  \qed

\vskip 10pt

\noindent {\bf Proof of Theorem \ref{tricyclic-thm}: }
  Since $K_4$ is the unique graph in $\mathcal{G}_{4,6}$, the theorem holds for $n=4$.
 By Table 1 of \cite{CDS}, there are four $(5,7)$-graphs. By simple calculation, we can prove the theorem for $n=5$.
 By Table 1 of \cite{CP}, there are 22 $(6,8)$-graphs. By direct computation, we can prove the theorem for $n=6$.
 Now suppose that $n\geq 7$. By Lemma \ref{tricyclic-minimal energy-1}, it suffices to prove that
 $E(G)>E(S_{n,n+2})$ when $n\geq 10$, and $E(G)>E(B_{n,n+2})$ when $7\leq n\leq 9$ for $G\in \mathcal{G}^2_{n,n+2}$.

 Suppose that $G\in \mathcal{G}^2_{n,n+2}$ and $C_p$, $C_q$ are two disjoint odd cycles with $p+q\equiv 2$ (mod $4)$. Then there are at most two edge disjoint paths in $G$
 connecting $C_p$ and $C_q$.

 {\bf Case 1.} There exists exactly an edge disjoint path $P$ connecting $C_p$ and $C_q$.
 Then there exists an edge $e$ of $P$ such that $G-e=G_1\cup G_2$,
 where $G_1$ is an non-bipartite bicyclic graph with $n_1\geq 4$ vertices and $G_2$ is an non-bipartite unicyclic graph with $n_2\geq 3$ vertices.
 By Lemmas \ref{edge-cut}, \ref{lemma in Zhang}, \ref{unicyclic-minimal energy-1}, \ref{bicyclic-n-5,6,7} and Theorem \ref{bicyclic-thm}, we have
\begin{equation*}
\begin{split}
E(G)
&\geq    E(G_1\cup G_2) \geq   E(S_{n_1,n_1+1}\cup S_{n_2,n_2}) >  E(S_{n_1,n_1}\cup S_{n_2,n_2 }) \\
&\geq    E(S_{n-3,n-3}\cup C_3)  >   E(S_{n,n+2}).
\end{split}
\end{equation*}
 In particular, $E(G)> E(S_{n,n+2})> E(B_{n,n+2})$ for $7\leq n\leq 9$.

 {\bf Case 2.} There exist exactly two edge disjoint paths $P^1$ and $P^2$ connecting $C_p$ and $C_q$.
 Then there exist two edges $e_1$ and $e_2$ such that $e_i$ is an edge of $P^i$ for $i=1,2$, and $G-\{e_1,e_2\}=G_3\cup G_4$,
 where $G_3$ and  $G_4$ are non-bipartite unicyclic graphs.  Let $|V(G_3)|=n_1$ and
 $|V(G_4)|=n_2$.
 Then by Lemmas \ref{edge-cut}, \ref{lemma in Zhang} and \ref{unicyclic-minimal energy-1}, we have
 $$E(G)\geq E(G_3\cup G_4)\geq E(S_{n_1,n_1}\cup S_{n_2,n_2 })\geq E(S_{n-3,n-3}\cup C_3) > E(S_{n,n+2}).$$
 In particular, $E(G)> E(S_{n,n+2})> E(B_{n,n+2})$ for $7\leq n\leq 9$.
 The proof is thus complete. \qed

 \begin{rem}\label{tricyclic-remark}
 The proof of Theorem \ref{tricyclic-thm} (for large $n$) is similar to that of Lemma \ref{tricyclic-minimal
 energy-2} except that in \cite{ZK}, the authors did not point out that $G_1$ and
 $G_2$ are non-bipartite graphs.
\end{rem}

\section{The graph with minimal energy in $\mathcal{G}_{n,n+3}$}

 Li and Li \cite{LL} discussed the graph with minimal energy in $\mathcal{G}^1_{n,n+3}$, and
 we first restate their results.

\begin{figure}[ht]
\centering
    \setlength{\unitlength}{0.05 mm}%
  \begin{picture}(3006.8, 546.9)(0,0)
  \put(0,0){\includegraphics{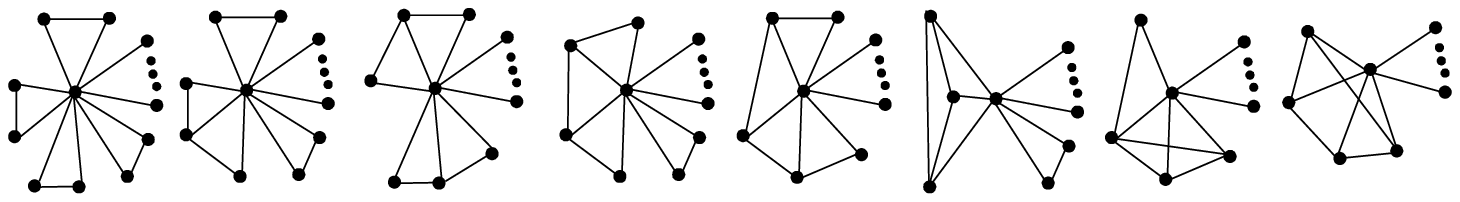}}
  \put(123.15,38.20){\fontsize{8.53}{10.24}\selectfont \makebox(150.0, 60.0)[l]{$G_1$\strut}}
  \put(477.36,40.13){\fontsize{8.53}{10.24}\selectfont \makebox(150.0, 60.0)[l]{$G_2$\strut}}
  \put(835.41,38.20){\fontsize{8.53}{10.24}\selectfont \makebox(150.0, 60.0)[l]{$G_3$\strut}}
  \put(1256.99,40.13){\fontsize{8.53}{10.24}\selectfont \makebox(150.0, 60.0)[l]{$G_4$\strut}}
  \put(1605.42,43.98){\fontsize{8.53}{10.24}\selectfont \makebox(150.0, 60.0)[l]{$G_5$\strut}}
  \put(1998.13,40.13){\fontsize{8.53}{10.24}\selectfont \makebox(150.0, 60.0)[l]{$G_6$\strut}}
  \put(2367.73,42.05){\fontsize{8.53}{10.24}\selectfont \makebox(150.0, 60.0)[l]{$G_7$\strut}}
  \put(2733.49,43.98){\fontsize{8.53}{10.24}\selectfont \makebox(150.0, 60.0)[l]{$G_8$\strut}}
  \end{picture}%
  \caption{Graphs $G_1, G_2, G_3, G_4, G_5, G_6, G_7$ and $G_8$.}\label{fig-tetracyclic}
\end{figure}

 Follow \cite{LL}, let $G_1, G_2, \ldots, G_8$ be
 eight special graphs in $\mathcal{G}_{n,n+3}$ as shown in Figure \ref{fig-tetracyclic}.
 Let $\mathscr{I}_n=\{S_{n,n+3}, B_{n,n+3}, G_1, G_2, G_3,G_4,G_5,G_6,G_7,G_8\}$.

 \begin{lem}\textnormal{\cite{LL}}\label{tetracyclic-minimal energy-1}
 If $G\in \mathcal{G}^1_{n,n+3}$ and $G\not\in \mathscr{I}_n$, then $E(G)>E(B_{n,n+3})$ for $n\geq 9$.
\end{lem}

  In fact, Lemma \ref{tetracyclic-minimal energy-1} is also true for
  $n=8$.

 \begin{lem}\label{tetracyclic-minimal energy-2}
 If $G\in \mathcal{G}^1_{8,11}$ and $G\not\in \mathscr{I}_8\setminus\{G_1\}$, then $E(G)>E(B_{8,11})$.
\end{lem}

 \pf By the results (see the proofs of Lemma 2.2 and Proposition 2.3) of \cite{LL},
 all we need is to show that $b_4(G)-b_4(B_{8,11})>0$ when $G$ contains exactly $i$ ($i=10,12,13,14,15$)
 cycles (see Case 7 of Lemma 2.2). From \cite{LL}, we have
 $$b_4(G)-b_4(B_{8,11})\geq \frac{1}{2}n^2+\frac{3}{2}n-12-2s-(5n-35),$$
 where $s$ is the number of quadrangles in $G$.
 It is easy to check that in this case, $G$ has at most $13$ quadrangles.
 %%%%Note that for any graph $G$, $b_4(G)=m(G,2)-2s$, where $m(G,2)$ is the number of $2$-matchings of $G$ and $s$ is the number of quadrangles in $G$ \cite{ZZ}.
 Therefore $$b_4(G)-b_4(B_{8,11})\geq \frac{1}{2}n^2+\frac{3}{2}n-12-26-(5n-35)=\frac{1}{2}n(n-7)-3=1>0.$$
 The proof is thus complete. \qed

 From Lemma \ref{lemma Sn,e and Bn,e}, we can obtain the following result.

\begin{cor}\label{tetracyclic-minimal energy-3}
 $E(S_{n,n+3}) < E(B_{n,n+3})$ for $n\geq 12$, and $E(B_{n,n+3}) < E(S_{n,n+3})$ for
 $7\leq n\leq 11$.
\end{cor}

 In \cite{LL}, the authors failed to get the above result in that (in the proof of Proposition 2.5 of \cite{LL})
  they used the wrong formula $b_4(S_{n,n+3})=4n-18$ instead of the
  correct one
  $b_4(S_{n,n+3})=4n-24$. They also gave the following result.

\begin{lem}\textnormal{\cite{LL}}\label{tetracyclic-minimal energy-4}
 For each $G_j\in \mathscr{I}_n$ ($j=1,\ldots,8$), $E(S_{n,n+3})
  <E(G_j)$ for $n\geq 9$ and $E(B_{n,n+3})<E(G_j)$ for $9\leq n\leq
  17$.
\end{lem}

 By the proof of Lemma 2.4 of \cite{LL}, we can get the following result
 for $n=8$.

 \begin{lem}\label{tetracyclic-minimal energy-5}
 For each $G_j\in \mathscr{I}_n\setminus\{G_1\}$ ($j=2,\ldots,8$), $E(B_{n,n+3})<E(G_j)$ for $n=8$.
\end{lem}

  By Lemmas \ref{tetracyclic-minimal energy-1}, \ref{tetracyclic-minimal energy-2}, \ref{tetracyclic-minimal energy-4},\ref{tetracyclic-minimal
  energy-5} and Corollary \ref{tetracyclic-minimal energy-3}, we can
  characterize the graph with minimal energy in $\mathcal {G}^1_{n,n+3}$.

\begin{lem}\label{tetracyclic-minimal energy-6}
 The graph with minimal energy in $\mathcal {G}^1_{n,n+3}$ is $B_{n,n+3}$ for $8\leq n\leq 11$, and $S_{n,n+3}$ for $n\geq 12$, respectively.
\end{lem}

  %%%To give a complete solution to Conjecture \ref{conjecture-minimal energy} for $e=n+3$, we need the following two lemmas.
 To prove Theorem \ref{tetracyclic-thm}, we need the following two lemmas.

 \begin{lem}\label{lemma compare tricyclic with unicyclic}
 (1) $E(K_4) > E(S_{4,4})$, and $E(B_{n,n+2}) > E(S_{n,n})$ for $7\leq n\leq 9$.

 (2) $E(S_{n,n+2}) > E(S_{n,n})$ for $n\geq 5$.
\end{lem}

 \pf (1) It is easy to obtain that
 $E(K_4)=6$, $E(S_{4,4})\doteq 4.96239$, $E(B_{7,9})\doteq 7.21110$, $E(S_{7,7})\doteq 6.64681$,
 $E(B_{8,10})\doteq 7.91375$, $E(S_{8,8})\doteq 7.07326$,
 $E(B_{9,11})\doteq 8.46834$ and $E(S_{9,9})\doteq 7.46410$. Hence
 the result (1) follows.

 (2) Since $6=E(S_{5,7}) > E(S_{5,5})\doteq 5.62721$, we now suppose
 $n\geq 6$. By direct computation, we have that
 $\phi(S_{n,n+2},x)=x^n-(n+2)x^{n-2}-6x^{n-3}+(3n-15)x^{n-4}$ and
 $\phi(S_{n,n},x)=x^n-nx^{n-2}-2x^{n-3}+(n-3)x^{n-4}$. By Eq.
 \eqref{energy-1}, we obtain that
 \begin{eqnarray*}
 E(S_{n,n+2}) & = &\frac{1}{2\pi}{\large\int}_{-\infty}^{+\infty}\frac{1}{x^{2}}\ln((1+(n+2)x^2+(3n-15)x^4)^2+(6x^3)^2) dx\\
              & > &\frac{1}{2\pi}{\large\int}_{-\infty}^{+\infty}\frac{1}{x^{2}}\ln((1+nx^2+(n-3)x^4)^2+(2x^3)^2) dx\\
              & = & E(S_{n,n}).
\end{eqnarray*}
\qed

\begin{lem}\label{lemma compare tetracyclic with unicyclic}
  $E(S_{n-3,n-3}\cup C_3) > E(S_{n,n+3})$ for $n\geq 6$.
\end{lem}

 \pf For $6\leq n\leq 14$, the result follows by direct computation.
 Suppose that $n\geq 15$.
 By direct calculation, we have that
$\phi(S_{n,n+3},x)=x^n-(n+3)x^{n-2}-8x^{n-3}+(4n-24)x^{n-4}$. Let
 $f(x)=x^4-(n+3)x^{2}-8x+4n-24$. Then we have that
 $f(-\sqrt{n-1})> 0$,  $f(-2)<0$, $f(0)>0$,  $f(2)<0$ and $f(\sqrt{n+3})> 0$.
 Hence $$E(S_{n,n+3})<4+\sqrt{n-1}+\sqrt{n+3}.$$
 On the other hand, we have $E(S_{n-3,n-3}\cup C_3)> 4+ \sqrt{2}+2\sqrt{n-4}$
 \cite{ZK}, and so $E(S_{n-3,n-3}\cup C_3) > E(S_{n,n+3})$.
 \qed

\vskip 10pt

\noindent {\bf Proof of Theorem \ref{tetracyclic-thm}: }
 By Table 1 of \cite{CDS}, there are two $(5,8)$-graphs. By simple calculation, we can prove the theorem for $n=5$.
 By Table 1 of \cite{CP}, there are 20 $(6,9)$-graphs. By direct computation, we can prove the theorem for $n=6$.
 By \cite{CDGT}, there are 132 $(7,10)$-graphs. By direct computing, we can prove the theorem for $n=7$.
 Now suppose that $n\geq 8$. By Lemma \ref{tetracyclic-minimal energy-6} and Corollary \ref{tetracyclic-minimal energy-3}, it suffices to prove that
 $E(G)>E(S_{n,n+3})$ for $G\in \mathcal{G}^2_{n,n+3}$.

 Suppose that $G\in \mathcal{G}^2_{n,n+3}$ and $C_p$, $C_q$ are two disjoint odd cycles with $p+q\equiv 2$ (mod $4)$. Then there are at most three edge disjoint paths in $G$
 connecting $C_p$ and $C_q$.

 {\bf Case 1.} There exists exactly an edge disjoint path $P^1$ connecting $C_p$ and $C_q$.
 Then there exists an edge $e_1$ of $P^1$ such that $G-e_1=G_1\cup G_2$,
 where either both $G_1$ and $G_2$ are non-bipartite bicyclic graphs, or $G_1$ is an non-bipartite tricyclic graph and $G_2$ is an non-bipartite unicyclic graph.
 Let $|V(G_1)|=n_1$ and $|V(G_2)|=n_2$.

  {\bf Subcase 1.1.} Both $G_1$ and $G_2$ are non-bipartite bicyclic graphs.
  Then by Lemmas \ref{edge-cut}, \ref{lemma in Zhang}, \ref{bicyclic-n-5,6,7}, \ref{lemma compare tetracyclic with unicyclic} and Theorem \ref{bicyclic-thm}, we have
\begin{equation*}
\begin{split}
E(G)
& \geq   E(G_1\cup G_2) \geq   E(S_{n_1, n_1+1}\cup S_{n_2, n_2+1 })  >  E(S_{n_1 ,n_1}\cup S_{n_2, n_2 })\\
& \geq    E(S_{n-3, n-3}\cup C_3) > E(S_{n, n+3}).
\end{split}
\end{equation*}

 {\bf Subcase 1.2.} $G_1$ is an non-bipartite tricyclic graph and $G_2$ is an non-bipartite unicyclic graph.
 It follows from Theorem \ref{tricyclic-thm} and Lemma \ref{lemma compare tricyclic with
 unicyclic} that $E(G_1)> E(S_{n_1,n_1})$. Therefore
  by Lemmas \ref{edge-cut}, \ref{lemma in Zhang}, \ref{unicyclic-minimal energy-1} and \ref{lemma compare tetracyclic with unicyclic}, we have
 $$E(G) \geq  E(G_1\cup G_2)  >  E(S_{n_1,n_1}\cup S_{n_2,n_2 }) \geq E(S_{n-3,n-3}\cup C_3) >
 E(S_{n,n+3}).$$

 {\bf Case 2.} There exist exactly two edge disjoint paths $P^2$ and $P^3$ connecting $C_p$ and $C_q$.
 Then there exist two edges $e_2$ and $e_3$ such that $e_i$ is an edge of $P^i$ for $i=2,3$, and $G-\{e_2,e_3\}=G_3\cup G_4$,
 where $G_3$ is an non-bipartite bicyclic graph with $n_1$ vertices and $G_4$ is an non-bipartite unicyclic graph with $n_2$ vertices.
 By Lemmas \ref{edge-cut}, \ref{lemma in Zhang}, \ref{unicyclic-minimal energy-1}, \ref{bicyclic-n-5,6,7}, \ref{lemma compare tetracyclic with unicyclic} and Theorem \ref{bicyclic-thm}, we have
\begin{equation*}
\begin{split}
E(G)
& \geq   E(G_3\cup G_4)\geq   E(S_{n_1,n_1+1}\cup S_{n_2,n_2 }) > E(S_{n_1,n_1}\cup S_{n_2,n_2})\\
& \geq    E(S_{n-3,n-3}\cup C_3) >E(S_{n,n+3}).
\end{split}
\end{equation*}

 {\bf Case 3.} There exist exactly three edge disjoint paths $P^4$, $P^5$ and $P^6$ connecting $C_p$ and $C_q$.
 Then there exist three edges $e_4$, $e_5$ and $e_6$ such that  $e_i$ is an edge of $P^i$ for $i=4,5,6$, and $G-\{e_4,e_5,e_6\}=G_5\cup G_6$,
 where $G_5$ and  $G_6$ are non-bipartite unicyclic graphs.
  Let $|V(G_5)|=n_1$ and $|V(G_6)|=n_2$.
 Then by Lemmas \ref{edge-cut}, \ref{lemma in Zhang}, \ref{unicyclic-minimal energy-1} and \ref{lemma compare tetracyclic with unicyclic}, we have
 $$E(G)\geq E(G_5\cup G_6)\geq E(S_{n_1,n_1}\cup S_{n_2,n_2 })\geq E(S_{n-3,n-3}\cup C_3) > E(S_{n,n+3}).$$

 The proof is thus complete. \qed

\vskip 0.5cm

\noindent{\bf Acknowledgments} \vskip 0.4cm

  This work is supported by  NNSFC (Nos. 11101351 and 11171288), NSF of the Jiangsu Higher Education Institutions (No. 11KJB110014) and the government scholarship of Jiangsu
  province.

\end{document}